\documentclass[12pt]{article}
\usepackage{graphicx}
\usepackage{subfigure}
 \usepackage{hangcaption}
\usepackage{cite}
 \usepackage{amsmath,amssymb}

\begin{document}
\def\roughly#1
{\raise.3ex\hbox{$#1$\kern-.75em\lower1ex\hbox{$\sim$}}}
\baselineskip0.32in

\vspace{-2cm}

 \def\half{ {1\over 2} }
 \def\tab{ {\hskip 0.15 true in} }
 \def\vtab{ {\vskip 0.1 true in} }
 \def\vtabb{ {\vskip 0.0 true in} }
 \def\blah{ {\vskip 0.1 true in} }
 \def\Order{\mbox{$\cal{O}$}}
 \def\vec#1{ {\overline {#1}} }
 \def\vecc#1{ {\overline {#1}} }
 \def\congruent{=}

      \newcommand{\beq}[1]{ \begin{equation} \label{eq.#1} }
      \newcommand{\eeq}{ \end{equation} }
       \newcommand{\IR}{I \!\! R}
       \newcommand{\IS}{I \!\!\! S}
       \newcommand{\IZ}{Z}
       \newcommand{\sfrac}[2]{{\scriptstyle \frac{#1}{#2}}}
        \newcommand{\RR}{\bf R}
        \newcommand{\XX}{\bf X}
\newcommand{\bfomega}{\mbox{\boldmath $\omega$ \unboldmath} \hskip -0.075 true in}
\newcommand{\xx}{{\bf x}}
\newcommand{\yy}{{\bf y}}
\newcommand{\bb}{{\bf b}}
\newcommand{\zz}{{\bf z}}
\newcommand{\qq}{{\bf q}}
\newcommand{\pp}{{\bf p}}
\newcommand{\bfxi}{\mbox{\boldmath $\xi$ \unboldmath} \hskip -0.075 true in}
\newcommand{\bfsigma}{\mbox{\boldmath $\sigma$ \unboldmath} \hskip -0.075 true in}
\newcommand{\bftau}{\mbox{\boldmath $\tau$ \unboldmath} \hskip -0.075 true in}
\newcommand{\bfepsilon}{\mbox{\boldmath $\epsilon$ \unboldmath} \hskip -0.075 true in}
\newcommand{\bfphi}{\mbox{\boldmath $\phi$ \unboldmath} \hskip -0.075 true in}
\newcommand{\bfchi}{\mbox{\boldmath $\chi$ \unboldmath} \hskip -0.075 true in}
\newcommand{\bfdelta}{\mbox{\boldmath $\delta$ \unboldmath} \hskip -0.075 true in}
\newcommand{\Real}{\ensuremath{\mathbb R}}
\newcommand{\SE}[1]{\ensuremath{{\mathrm{SE}(#1)}}}
\newcommand{\SO}[1]{\ensuremath{{\mathrm{SO}(#1)}}}

\vspace{0.5cm}

\centerline{ \large  {\bf Deblurring of Motionally Averaged Images
with }} \centerline{ \large  {\bf Applications to Single-Particle
Cryo-Electron Microscopy}}

\vspace{0.5cm}

\centerline{Wooram Park, Daniel N. Rockmore\footnote{Departments
of Mathematics and Computer Science, Dartmouth College},
% Peter T. Kim \footnote{Department of Mathematics and Statistics, University
% of Guelph, Ontario, Canada}
Dean Madden\footnote{Department of Biochemistry, Dartmouth
College} and Gregory S. Chirikjian}

\vspace{0.5cm}

\centerline{\it Department of Mechanical Engineering}
\centerline{\it  Johns Hopkins University} \centerline{\it
Baltimore, MD \ \ 21218, \ USA}

\vskip 0.05 true in

\centerline{\bf Abstract} \vskip 0.05 true in
%\begin{abstract}

This paper addresses the deconvolution of an image that has been
obtained by superimposing many copies of an underlying unknown
image of interest. The superposition is assumed to not be exact
due to noise, and is described using an error distribution in
position, orientation, or both. We assume that a good estimate of
the error distribution is known. The most natural approach to take
for the purely translational case is to apply the Fourier
transform and use the classical convolution theorem together with
a Weiner filter to invert. In the case of purely rotational
deblurring, the similar Fourier analysis is applied. That is, for
an blurred image function defined on polar coordinates, the
Fourier series and the convolution theorem for the series can be
applied. In the case of combinations of translational and
rotational errors, the motion-group Fourier transform is used. In
addition, for the three cases we present the alternative method
using Hermite and Laguerre-Fourier expansion, which has a special
property in Fourier transform. The problem that is solved here is
motivated by one of the steps in the cryo-electron-tomographic
reconstruction of biomolecular complexes such as viruses and ion
channels.
% \end{abstract}

\vskip 0.1 true in


\vskip 0.1 true in
%\noindent
%{\bf PACs Numbers:}

 \section{Introduction}

In single particle Cryo-Electron-Microscopy (Cryo-EM), the goal is
to reconstruct the 3D shape of large biomolecular complexes from
projection data. In particular, many essentially identical copies
of a complex of interest are embedded in a thin layer of vitreous
ice at randomized (and unknown) orientations. If we consider this
thin film of ice to be in the x-y plane in the lab frame, then an
electron beam takes projections of the density of the embedded
biomolecular complexes along the z direction. The goal is to
reconstruct the three-dimensional density of the complex (which
defines its shape) from these projection images. % \cite{baker}.
The difference between this problem and medical image
reconstruction is that the projection directions are unknown a
priori.

The signal-to-noise ratio in such measurements can be quite high
\cite{jfrankbook}. Therefore projections corresponding to the same
(or quite similar) projection directions are grouped together and
superimposed. In doing so, the random noise of the background has
a tendency to cancel, and the features of interest in the
projections reinforce each other as the number of superimposed
projections becomes large \cite{jfrankbook}. This averaging
technique may not be used in all types of 3D reconstruction of the
electron microscopy. However, it is still useful for the analysis
of classes of 2D projected images, the raw data of which contains
a large amount of noise.

One problem is that the superposition of images might not be
exact. This results in a blurring relative to the true underlying
image of interest. To get a sense of this, let us consider the
following image model including zero-mean ergodic white
noise\cite{jfrankbook}.
\begin{equation}
\rho({\bf x},t) = \rho_0( g_{t}^{-1}   {\bf x}) + n({\bf x},t).
\label{mainmodel}
\end{equation}
Here $g_t$ is the homogeneous transformation in $\SE{2}$ ( the Lie
group describing the translational and rotational motion),  $n$ is
the noise, ${\bf x} \in \Real^2$ is the planar position of points
in each image, and $ t \in \Real^{+}$ is an artificial time
variable used to order the images. If there is no noise term, it
is intuitively clear that the appropriate matching of two images,
$\rho({\bf x},t_i)$ and $\rho({\bf x},t_j)$ occurs at
$g_{t_i}=g_{t_j}$ and the superposition is $( \rho({\bf x},t_i) +
\rho({\bf x},t_j) )/2 = f_0( g_{t_i}^{-1}   {\bf x})$. However, if
we have the noise term, the matching of many data images produce
various $g_t$. The superposition of them will be in the form of
$$
\frac{1}{N} \sum^{N}_{i=1} \rho({\bf x},t_i) \approx \frac{1}{N}
\sum^{N}_{i=1} \rho_0( g_{t_i}^{-1}   {\bf x}) = \int_G
\left(\frac{1}{N} \sum^{N}_{i=1} \delta(g_{t_i}^{-1} \circ g)
\right) \rho_0( g^{-1} {\bf x}).
$$
The first equality assumes that the noise term is approximately
canceled out during the superimposition, and the second equality
shows that this superposition is a convolution on the group of
rigid-body motions of the plane, $G=\SE2$. The right hand side is
essentially the blurred version of $ \rho_0(    {\bf x})$
depending on the distribution of $g_{t_i}$.

We therefore seek to solve the following inverse problem: Given a
blurred image, $\gamma({\bf x})$, that describes the optimal
superimposition of many experimentally-obtained projection images,
and given an estimate of the probability density function
describing the error distribution in the alignment of these
superimposed images, $f(g)$, we seek to find the deblurred image
$\rho({\bf x})$. This is expressed as the solution to the problem:
\begin{equation}
\int_{G} f(g) \rho(g^{-1} \cdot {\bf x}) dg = \gamma({\bf x}).
\label{main}
\end{equation}
Here $G$ is the group of transformations involved in alignment, $g
\cdot x$ denotes the group action of $G$ on $\Real^2$, and $dg$ is
the associated invariant integration measure for that group
\cite{book}. In this paper we consider three cases: (1) $G =
(\Real^2,+)$, the translation group in the plane; (2) $G=SO(2)$,
the rotation group in the plane; and (3) $G=\SE2$, the Euclidean
motion group of the plane. Explicitly, in these three cases we
have

\begin{equation}
\int_{-\infty}^{\infty} \int_{-\infty}^{\infty} f_1(y_1, y_2)
\rho(x_1 - y_1, x_2 - y_2) dy_1 dy_2 = \gamma_1(x_1, x_2),
\label{main1}
\end{equation}

\begin{equation}
 \int_{0}^{2\pi} f_2(\theta) \rho(x_1 \cos \theta + x_2 \sin \theta; -x_1 \sin \theta + x_2 \cos \theta )
d\theta = \gamma_2(x_1, x_2), \label{main2}
\end{equation}

\begin{equation}
\int_{-\infty}^{\infty}  \int_{-\infty}^{\infty} \int_{0}^{2\pi}
f_3(y_1,y_2,\theta) \rho((x_1-y_1) \cos \theta + (x_2 - y_2)\sin
\theta; -(x_1-y_1) \sin \theta + (x_2 - y_2) \cos \theta ) dy_1
dy_2 d\theta  \label{main3}
\end{equation}
$$
= \gamma_3(x_1, x_2).
$$

As a model for the functions $f_i$, we will assume appropriate
concepts of Gaussian distributions. Recall that the diffusion
equation on the line,
$$ \frac{\partial f}{\partial t} = \frac{\partial^2 f}{\partial x^2}, $$
subject to the initial conditions $f(x,0) = \delta(x)$, has the
solution
\begin{equation}
f(x,t) = \frac{1}{2\sqrt{\pi t}} e^{-x^2/4t}. \label{gauss1}
\end{equation}
The solution to the uniform planar diffusion equation:
$$ \frac{\partial f_1}{\partial t} = \frac{\partial^2 f_1}{\partial x_1^2}
+ \frac{\partial^2 f_1}{\partial x_1^2}, $$ can be written as
$$ f_1(x_1, x_2, t) = f(x_1,t) f(x_2,t). $$
Likewise, the diffusion on the circle can be viewed as a folded
normal distribution:
$$ f_2(\theta,t) = \sum_{n=-\infty}^{\infty} f(\theta - 2\pi n, t). $$
It is useful to note that this can be expressed alternatively as a
Fourier series:
$$ f_2(\theta,t) = \frac{1}{2\pi} \sum_{k=-\infty}^{\infty} e^{-k^2 t} e^{i k \theta}. $$
A model for combined translational and rotational error is then
$$ f_3(x_1,x_2,\theta; t_1,t_2) = f_1(x_1, x_2, t_1) f_2(\theta,t_2) $$
where the small values of $t_1$ and $t_2$ can be chosen separately
to describe different amounts of translational and rotational
error, as well as to account for the fact that the units of
measurement are different for translations and rotations.

\section{Deconvolution of Motion-Averaged Images Using Fourier Transform}

%%%%%%%% 2 %%%%%%%%%%%%%%%

In this section we address how to solve each of the three
deconvolution problems in (\ref{main1}) - (\ref{main3}) using
Fourier transform.

\subsection{The 2-D Fourier Transform for Translational Deconvolution}
\label{sc:subsc_FT_case_1}

%%%%%%%% 2.1 %%%%%%%%%%%%%%%

The natural tool to use to solve the deconvolution problem in
(\ref{main1}) is the Fourier transform. The Fourier transform in
two dimensions is written in Cartesian coordinates as
$$ \hat{f}({\bf \omega})  = \hat{f}(\omega_1,\omega_2) =
\int_{-\infty}^{\infty} \int_{-\infty}^{\infty} f(x_1,x_2)
e^{-i(\omega_1 x_1 + \omega_2 x_2)} dx_1 dx_2 = \int_{\Real^2}
f({\bf x}) e^{-i {\bf \omega} \cdot {\bf x}} d{\bf x} $$ and the
inversion formula is
$$ f({\bf x}) = \frac{1}{(2\pi)^2}
\int_{\Real^2} \hat{f}({\bf \omega}) e^{i {\bf \omega} \cdot {\bf
x}} d{\bfomega}
$$
Fourier transform of the distribution function, $f_1(x_1,x_2,t)$
in (\ref{main1}) is $ \hat{f}_1({\bf \omega}) =
e^{-(\omega_1^2+\omega_2^2)t} $.

The convolution theorem then converts (\ref{main1}) to a problem
in Fourier space of the form
$$ \hat{f}_1({\bf \omega}) \hat{\rho}({\bf \omega}) = \hat{\gamma}_1({\bf \omega}), $$
which is inverted after regularization as
\begin{equation}
 \hat{\rho}({\bf \omega}) = \hat{\gamma}_1({\bf \omega}) \overline{\hat{f}_1({\bf \omega})}
 /(\epsilon + |\hat{f}_1({\bf \omega})|^2).
\label{fourprod}
\end{equation}
The regularization parameter, $\epsilon$, is a very small positive
number that is introduced to handle zeros of the Fourier
transform. In fact, for the Gaussian distribution of interest in
our problem, there are no zeros, but in both real and Fourier
space the tails of the distribution can approach zero at points
sufficiently far from the origin. (\ref{fourprod}) is nothing more
than the well-known Wiener filter.

\subsection{Deconvolution of Purely Rotational Misalignment}
\label{sc:subsc_FT_case_2}
%%%%%%%% 2.2 %%%%%%%%%%%%%%%

If the image functions are defined on polar coordinates,
(\ref{main2}) can be rewritten as
$$
 \int_{0}^{2\pi} f_2(\theta) \rho(r,\phi-\theta) d\theta =
 \gamma_2(r,\phi),
$$
where $x_1=r\cos\phi$ and $x_2=r\sin\phi$. If we fixed the value
of $r$, (\ref{main2}) becomes the convolution of the two functions
on a circle as

\begin{equation}
 \int_{0}^{2\pi} f_2(\theta)   \rho^{(r)}( \phi-\theta) d\theta =
  \gamma^{(r)}_2(\phi),\label{eq:conv_so2}
\end{equation}

where $ \rho^{(r)}( \phi-\theta) = \rho(r,\phi-\theta) $ and $
\gamma^{(r)}_2(\phi) = \gamma_2(r,\phi).$

The Fourier series expansion of a function defined on a circle
gives
$$
f(\theta)= \frac{1}{2\pi} \sum^{ \infty }_{-\infty} f_n e^{ i n
\theta },
$$
where
$$
f_n =  \int^{2\pi}_{0}f(\theta) e^{-i n t } d \theta.
$$
The Fourier transform of the distribution function,
$f_2(\theta,t)$ in (\ref{main2}) is $ (\hat{f}_2)_n = e^{-n^2 t}
$.

The convolution theorem of Fourier series converts
(\ref{eq:conv_so2}) to the problem in Fourier space of the form
$$
(\hat{f}_2)_n  (\hat{\rho}^{(r)})_n  = (\hat{\gamma}^{(r)})_n.
$$
As in the case of the translational deconvolution, the inversion
with regularization is
\begin{equation}\label{eq:ft_case_2_reg}
 (\hat{\rho}^{(r)})_n  =  (\hat{\gamma}^{(r)})_n   \overline{(\hat{f}_2)_n}
 /(   \epsilon  + |     (\hat{f}_2)_n     | ^2   )     .
\end{equation}

\subsection{Deconvolution of Combined Translational and Rotational Blurring}
%%%%%%%% 2.3 %%%%%%%%%%%%%%%
\label{sc:subsc_FT_case_3}

 In order to solve the full motional
deconvolution problem, the appropriate concept of Fourier
transform is required. In particular, since $f_3$ is a function on
the group of rigid-body motions of the plane, $\SE2$, and a
function on $\Real^2$ can be viewed as a function on $\SE2$ that
is constant over the orientational variable, (\ref{main3}) can be
viewed as a convolution on $\SE2$. We therefore review here the
group $\SE2$ and the associated Fourier analysis.

\subsubsection{Representation Theory of The Euclidean Motion Group of the Plane}
%%%%%%%% 2.3.1 %%%%%%%%%%%%%%%

Each element of $\SE2$ is parameterized in either rectangular or
polar coordinates as:
$$ g(a_1, a_2,\theta) = \left(\begin{array}{ccc}
\cos \theta & -\sin \theta & a_1 \\
\sin \theta & \cos \theta & a_2 \\
0 & 0 & 1
\end{array}\right) $$
or
$$
g(a,\phi,\theta) = \left(\begin{array}{ccc}
\cos \theta & -\sin \theta & a\cos \phi \\
\sin \theta & \cos \theta & a\sin \phi \\
0 & 0 & 1
\end{array}\right), $$
where $a =\| {\bf a} \|$.

A irreducible unitary representations of $\SE2$ (see \cite{book,
10Vilenkin,10sugiura} for general definition) can be viewed as
infinite dimensional matrices, $U(g,p)$ with elements expressed as
\begin{equation}
u_{mn}(g(a,\phi,\theta),p) = i^{n-m} e^{-i[n\theta + (m-n)\phi]}
J_{n-m}(p \, a) \label{elements1}
\end{equation}
where $J_{\nu}(x)$ is the $\nu^{th}$ order Bessel function and $m$
and $n$ range over all integer values.

From this expression, and the fact that ${U}(g,p)$ is a unitary
representation, we have that:
$$ u_{mn}(g^{-1}(a,\phi,\theta),p) = u_{mn}^{-1}(g(a,\phi,\theta),p) = $$
\begin{equation}
\overline{u_{nm}(g(a,\phi,\theta),p)} = i^{n-m} e^{i[m\theta +
(n-m)\phi]} J_{m-n}(p a). \label{se2elem2}
\end{equation}

\vspace{0.5cm}

These matrix elements are related by the symmetries:
\begin{equation}
\label{sym2d} \overline{u_{m n}(g,p)}\,=\, (-1)^{m-n}
u_{-m,-n}(g,p),
\end{equation}
\begin{equation}
\label{symtwopp} {u_{m
n}(g(-a,\phi,\theta),p)}\,\stackrel{\triangle}{=} \, {u_{m
n}(g(a,\phi \pm \pi,\theta),p)}\,=\, (-1)^{m-n}
u_{m,n}(g(a,\phi,\theta),p)
\end{equation}
and
\begin{equation}
\label{unitarpp} (-1)^{m-n} u_{m,n}(g(a,\phi -\theta,-\theta),p) =
\overline{u_{n m}(g(a,\phi,\theta),p)}.
\end{equation}
The equality in (\ref{unitarpp}) follows from (\ref{se2elem2}) and
(\ref{symtwopp}).

\subsubsection{The Fourier Transform for the Euclidean Motion Group of the Plane}
%%%%%%%% 2.3.2 %%%%%%%%%%%%%%%

The Fourier transform of a sufficiently well-behaved function on
$\SE2$, and the corresponding inverse transform are defined as:
$$ {\cal F}(f) = \hat{f}(p) = \int_{G} f(g) U(g^{-1},p) \, dg $$
and
$$ {\cal F}^{-1}(\hat{f}) =
f(g) = \int_{\,0}^{\,\infty} {\rm trace}(\hat{f}(p) U(g,p)) p dp.
$$

As with the Fourier transform of functions on ${\Real}^N$,
$$ {\cal F} {\cal F}^{-1}(\hat{f}) = \hat{f} \tab \tab {\cal F}^{-1} {\cal F} (f) = f. $$
A proof that these identities hold is given in \cite{10sugiura}.
The fact that the inverse transform works depends on $\{U(g,p)\}$
being a complete set of irreducible representations, and the fact
that it is unitary allows us to write $U(g^{-1},p) =
U^{\dagger}(g,p)$ instead of computing the inverse of an infinite
dimensional matrix.

The matrix elements of the transform can be calculated using the
matrix elements of $U(g,p)$ defined in (\ref{elements1}) as:
\begin{equation}\label{eq:ft_se2}
\hat{f}_{mn}(p) = \int_{G} f(g) u_{mn}(g^{-1},p) \, dg.
\end{equation}
Likewise, the inverse transform can be written in terms of
elements as:
$$ f(g) = \sum_{n,m \in \IZ}
\int_{0}^{\infty} \hat{f}_{mn}(p) u_{nm}(g,p) p dp. $$

\subsubsection{Regularized Deconvolution of Motional Deblurring in $\SE2$}
%%%%%%%% 2.3.3 %%%%%%%%%%%%%%%

Given motional blurring expressed in (\ref{main}) when $G=\SE2$,
the result becomes (\ref{main3}). This is a convolution on the
motion group. The result can be solved by applying the motion
group Fourier transform to yield:
\begin{equation}\label{eq:ft_conv_3}
 \hat{\rho}(p) \hat{f}_3(p) = \hat{\gamma}_3(p).
 \end{equation}
The functions $\hat{\rho}(p)$ and $\hat{\gamma}_3(p)$ are row
vectors. The direct matrix inversion would be
\begin{equation}\label{eq:ft_case_3_reg}
\hat{\rho}(p) = \hat{\gamma}_3(p) [\hat{f}_3(p)]^{-1} .
\end{equation}
However, if the matrix $\hat{f}_3(p)$ becomes singular, then this
needs to be regularized. The procedure for doing this is explained
in  \cite{fredholm} , and involves the computation of a weighted
least-squares pseudo-inverse.

In the current context, we can compute the entries of
$\hat{f}_3(p)$ analytically. For the time being, we drop the
subscript `3', and write the function in polar coordinates as
$$
 f(r,\phi,\theta; t_1,t_2) = \frac{1}{8{\pi^2 t_1}} e^{-r^2/4t_1} \sum_{k=-\infty}^{\infty}
e^{-k^2 t_2} e^{i k \theta}.
$$
We use the fact that in polar coordinates $dg = rdr d\phi d\theta$
and the above function is independent of $\phi$ (which will result
in a diagonal $\SE2$ Fourier transform matrix). Computing the
$\SE2$ Fourier transform of the distribution, we find:
\begin{equation} \label{eq:ft_f3}
 \hat{f}_{mn}(p) = \frac{1}{2 t_1} \delta_{mn} \left(\int_{0}^{\infty} e^{-r^2/4t_1}
 J_0(pr) r dr\right) e^{-m^2 t_2}  =  \delta_{mn} e^{-p^2 t_1} e^{-m^2
 t_2}
\end{equation}
Therefore, the matrix, $\hat{f}_3(p)$ in (\ref{eq:ft_case_3_reg})
is diagonal. Its inversion is the simple inversion of scalar
values and the inversion with the regularization parameter is the
same as that in the previous two cases.

\section{Deconvolution of Motion-Averaged Images Using Hermite and Laguerre functions}
%%%%%%%% 3 %%%%%%%%%%%%%%%
Even though the Fourier transform is a good way to solve the
deconvolution problem, its implementation needs several
manipulations of data such as interpolation. We will discuss the
details in the next section. In this section, we develop the
alternative method of deconvolution using Hermite and Laguerre
functions. We utilize special properties that Hermite and
Laguerre-Fourier expansions have in the Fourier transform. In this
method, (\ref{fig:deblurred_images_case_1}) and
(\ref{fig:deblurred_images_case_3}) are solved in the Fourier
space, while (\ref{fig:deblurred_images_case_2}) is solved in the
real space.

Hermite function, $h_n(x)$ is an eigenfunction for the Fourier
transform as
$$
\int^{\infty}_{-\infty} h_n(x) e^{-i \omega x} dx = \sqrt{2 \pi}
(-i)^n h_n(\omega),
$$
where $h_n(x)$ is the Hermite function defined as
$$
h_n(x) = \frac{1}{s_n } H_n(x) e^{-x^2/2},
$$
where $s_n=\sqrt{2^n n! \sqrt{\pi}}$ and $H_n(x)$ is Hermite
polynomial, which is generated by the Rodrigues formula
$$
H_n(x) = (-1)^n e^{x^2} \frac{d^n}{dx^n}  (  e^{-x^2}    ).
$$
This property gives straightforward analytic solution, when the
image is defined as Hermite expansion and the Fourier method in
the previous section is applied.

While Hermite expansion is good for a function defined on
Cartesian coordinates, Laguerre-Fourier expansion looks better for
a function defined on polar coordinates. Mathematically, the two
expansions can be converted to each
other~\cite{park:interconversion}.

The {\it associated Laguerre polynomials} are generated by the
Rodrigues formula,
$$
L_n^k(x)=\frac{e^x x^{-k}}{n!}\frac{d^n}{dx^n}(e^{-x}x^{n+k}).
$$
The associated Laguerre polynomials are orthogonal over
$[0,\infty)$ with respect to the weighting function $x^k e^{-x}$,
$$
\int^{\infty}_{0} x^k e^{-x} L^{k}_{m}(x) L^{k}_{n}(x) dx =
\frac{(n+k)!}{n!}\delta_{m,n}.
$$
Using the Laguerre polynomials and Fourier basis, we can define
the basis function on two-dimensional polar coordinate as
follows~\cite{park:interconversion}~\cite{pola_rich}.
$$
\chi_{m,n}(r,\phi) = (-1) ^
{(m-|n|)/2}\sqrt{\frac{[(m-|n|)/2]!}{\pi[(m+|n|)/2]!}}r^{|n|}
L^{|n|}_{(m-|n|)/2}(r^2)e^{-r^2/2}e^{-in\phi},
$$
where $(m-|n|)$ and $(m+|n|)$ are even numbers. For convenience,
we sometimes divide it into two parts as
$$
\chi_{m,n}(r,\phi) = y_{m,n}(r)  z_n(\phi),
$$
where
\begin{eqnarray*}
y_{m,n}(r) &=& (-1) ^
{(m-|n|)/2}\sqrt{\frac{[(m-|n|)/2]!}{\pi[(m+|n|)/2]!}}r^{|n|}
L^{|n|}_{(m-|n|)/2}(r^2)e^{-r^2/2}, \\
z_n(\phi) &=& e^{-in\phi}
\end{eqnarray*}

 In this section we will show how to solve the aforementioned
deconvolution problem using the Hermite and Laguerre functions.
The appropriate expansion will be chosen for the three cases.

\subsection{The Translational Deconvolution using Hermite expansion}
\label{sc:subsc_HL_case_1}
%%%%%%%% 3.1 %%%%%%%%%%%%%%%

Hermite expansion of an image function is
$$
\rho(x_1,x_2) = \sum^{\infty}_{m=0} \sum^{\infty}_{n=0}
\check{\rho}_{mn}h_m(x_1) h_n(x_2),
$$
where
$$
\check{\rho}_{mn} = \int_{\Real^2} \rho(x_1,x_2)h_m(x_1) h_n(x_2)
dx_1 dx_2.
$$
When an image function can be expressed as a truncated Hermite
expansion with large $N$, the function is written as
\begin{equation}\label{eq:truncated_H}
\rho(x_1,x_2) =  \sum^{N}_{m=0} \sum^{N-m}_{n=0}
\check{\rho}_{mn}h_m(x_1) h_n(x_2).
\end{equation}
The Fourier transform of the image function is
$$
\hat{\rho}(\omega_1,\omega_2) = \sum^{N}_{m=0} \sum^{N-m}_{n=0}
\check{\rho}_{mn} (2\pi) (-i)^{m+n}  h_m(\omega_1) h_n(\omega_2).
$$
Since the Fourier transform of $f_1(x_1,x_2,t)$ is $
\hat{f}_1(\omega_1,\omega_2) = e^{-(\omega_1^2+\omega_2^2)t}$, the
convolution theorem gives
\begin{equation}
\hat{\gamma}_1(\omega_1,\omega_2) = \sum^{N}_{m=0}
\sum^{N-m}_{n=0} \check{\rho}_{mn} (2\pi) (-i)^{m+n} h_m(\omega_1)
h_n(\omega_2) e^{-(\omega_1^2+\omega_2^2)t}  .
\label{eq:HL_case1_LHS}
\end{equation}

On the other hand, since $H_m(x)$ is an $m$'th order polynomial,
it can be rewritten as
$$
H_m(x) = H_m(\frac{ax}{a}) = \sum^{m}_{k=0} \alpha_{m,k}(a^{-1})
H_k(ax),
$$
where $\alpha_{m,k}(a^{-1})$ is an appropriate coefficient
relating Hermite polynomial and its scaled version. Using this
expression, we can have
$$
h_m(\omega) e^{-\omega^2 t} = \sum^{m}_{k=0} \alpha_{m,k}(a^{-1})
\frac{s_k}{s_m} h_k(a\omega),
$$
where $a=\sqrt{2t+1}$. Therefore,
$$
\hat{\gamma}_1(\omega_1,\omega_2) = \sum^{N}_{m=0}
\sum^{N-m}_{n=0} \check{\rho}_{mn} (2\pi) (-i)^{m+n}
\sum^{m}_{k=0}\sum^{n}_{l=0}
\alpha_{m,k}(a^{-1})\alpha_{n,l}(a^{-1})
\frac{s_k}{s_m}\frac{s_l}{s_n} h_k(a\omega_1)h_l(a\omega_2) .
$$
We can reorder the summations and have
\begin{equation}
\hat{\gamma}_1(\omega_1,\omega_2) = \sum^{N}_{k=0}
\sum^{N-k}_{l=0} \left(  \sum^{N-l}_{m=k} \sum^{N-m}_{n=l}
\check{\rho}_{mn} (2\pi) (-i)^{m+n}
\alpha_{m,k}(a^{-1})\alpha_{n,l}(a^{-1})
\frac{s_k}{s_m}\frac{s_l}{s_n} \right)h_k(a\omega_1)h_l(a\omega_2)
. \label{eq:HL_case1_RHS}
\end{equation}
Its inverse Fourier transform is
$$
 \gamma_1 (x_1,x_2) = \sum^{N}_{k=0} \sum^{N-k}_{l=0}
\left(  \sum^{N-l}_{m=k} \sum^{N-m}_{n=l} \check{\rho}_{mn} (2\pi)
(-i)^{m+n} \alpha_{m,k}(a^{-1})\alpha_{n,l}(a^{-1})
\frac{s_k}{s_m}\frac{s_l}{s_n} \right)\frac{i^{k+l}}{2\pi
a^2}h_k(x_1/a)h_l(x_2/a) .
$$
Conversely, if we can have a truncated Hermite expansion for a
blurred image as
\begin{equation} \label{eq:blurred_H_case1}
 \gamma_1(x_1,x_2) =  \sum^{N}_{k=0} \sum^{N-k}_{l=0}
 \check{\gamma}_{kl} h_k(x_1/a)h_l(x_2/a),
\end{equation}
 then its Fourier transform is
\begin{equation}
\hat{\gamma}_1(\omega_1,\omega_2) = \sum^{N}_{k=0}
\sum^{N-k}_{l=0}
 \check{\gamma}_{kl} \frac{2\pi
a^2}{i^{k+l}}h_k(a\omega_1)h_l(a\omega_2) .
\label{eq:HL_case1_RHS_2}
\end{equation}
Equating (\ref{eq:HL_case1_LHS}) and (\ref{eq:HL_case1_RHS_2}) on
a various samples on $(\omega^{(p)},\omega^{(q)})$ gives
$$
(EHU)R(EHU)^T = H_a \times G \times H_a^T,
$$
where $E_{m,n}=\delta_{m,n} e^{-t(\omega^{(m)})^2}$,
$H_{m,n}=h_{n-1}(\omega^{(m)})$,
$(H_a)_{m,n}=h_{n-1}(a\omega^{(m)})$, $U_{m,n} = \delta_{m,n}
(-i)^m $, $R_{m,n}= \check{\rho}_{m-1,n-1}$ and $G_{m,n}=
\check{\gamma}_{m-1,n-1}$. In order to get $R$, which is the
Hermite coefficients for the deblurred image, we should examine
the inversion of the matrices.

Inversion of $U$ is given by $U^{-1}_{m,n} =\delta_{m,n} i^m
 $. Pseudo-inverse of $H$ is given by $H^{+}= (H^T
H)^{-1}H^T$ if the sampling points, $\omega^{(p)}$ are chosen
appropriately as shown in the previous
work~\cite{park:interconversion}. Inverting $E$ needs
regularization because inverse of $e^{-t(\omega^{(m)})^2}$ may be
unstable with large value of $\omega^{(m)}$. Therefore,
$E^{+}_{m,n}=\delta_{m,n} ( 1 / ( e^{-t(\omega^{(m)})^2} +
\epsilon )$ with a small number, $\epsilon$.

Now we have
$$
R = U^{-1} H^{+} E^{+} H_a \times G \times (U^{-1} H^{+} E^{+} H_a
)^T.
$$

\subsection{The Rotational Deconvolution using Laguerre-Fourier expansion}
\label{sc:subsc_HL_case_2}
%%%%%%%% 3.2 %%%%%%%%%%%%%%%

A 2D function defined on polar coordinate can be expressed as
$$
\rho(r,\phi) = \sum^{\infty}_{m=0}\sum^{m}_{n=-m}\tilde{\rho}_{mn}
\chi^*_{mn}(r,\phi),
$$
where
$$
\tilde{\rho}_{mn} = \int_{\Real^2} \rho(r,\phi) \chi_{mn}(r,\phi)
rdrd\phi
$$
If the function can be expressed as a truncated Laguerre-Fourier
expansion with large $N$, we have
\begin{equation} \label{eq:truncated_L}
\rho(r,\phi) = \sum^{N}_{m=0}\sum^{m}_{n=-m}\tilde{\rho}_{mn}
\chi^*_{mn}(r,\phi).
\end{equation}
Note that the integer variable $n$ increases by multiples of 2.

If the image function in (\ref{main2}) is defined on polar
coordinates, the convolution can be rewritten as
$$
\int_{0}^{2\pi} f_2(\theta) \rho(r,\phi-\theta ) d\theta =
\gamma_2(r,\phi),
$$
with the coordinate conversion, $x_1=r\cos\theta$ and
$x_2=r\sin\theta$. If we use the truncated Laguerre-Fourier
expansion for $\rho$, we have
$$
\int_{0}^{2\pi} \frac{1}{2 \pi}\sum^{\infty}_{k=-\infty} e^{-k^2
t} e^{ik\theta} \sum^{N}_{m=0}\sum^{m}_{n=-m}\tilde{\rho}_{mn}
\chi^*_{mn}(r,\phi-\theta )d\theta = \gamma_2(r,\phi),
$$
The left hand side can be computed as
$$
\frac{1}{2 \pi}
\sum^{\infty}_{k=-\infty}\sum^{N}_{m=0}\sum^{m}_{n=-m} e^{-k^2
t}\tilde{\rho}_{mn} y_{mn}(r) \int_{0}^{2\pi} e^{ik\theta}
e^{in(\phi-\theta)} d\theta =
\sum^{\infty}_{k=-\infty}\sum^{N}_{m=0}\sum^{m}_{n=-m} e^{-k^2 t}
e^{in\phi } \tilde{\rho}_{mn} y_{mn}(r)   \delta_{k,n}
$$
$$
= \sum^{N}_{m=0}\sum^{m}_{n=-m}  \left(\tilde{\rho}_{mn} e^{-n^2
t}\right) y_{mn}(r) z^*_{n}(\phi)= \sum^{N}_{m=0}\sum^{m}_{n=-m}
\left(\tilde{\rho}_{mn} e^{-n^2 t}\right) \chi^*_{mn}(r,\phi )
$$
Therefore the blurred image is
\begin{equation} \label{eq:blurred_L_case2}
\gamma_2(r,\phi)=\sum^{N}_{m=0}\sum^{m}_{n=-m} \tilde{\gamma}_{mn}
\chi^*_{mn}(r,\phi ),
\end{equation}
where
$$
\tilde{\gamma}_{mn} = \tilde{\rho}_{mn} e^{-n^2 t}.
$$
 This means that the convolved(blurred) image of a truncated
Laguerre-Fourier expansion with purely rotational motion is also a
truncated Laguerre-Fourier expansion with the same truncation
limit. The only difference is that the coefficients are scaled.
Once we compute the Laguerre-Fourier
coefficients($\tilde{\gamma}_{mn}$) of the blurred image, the
Laguerre-Fourier coefficients of the deblurred image is given by
\begin{equation}\label{eq:hl_case_2_reg}
\tilde{\rho}_{mn}= \tilde{\gamma}_{mn} \{ 1 / ( e^{-n^2 t}
+\epsilon )\}
\end{equation}
with regularization.

\subsection{The Translational and Rotational Deconvolution
using Laguerre-Fourier expansion} \label{sc:subsc_HL_case_3}
%%%%%%%% 3.3 %%%%%%%%%%%%%%%

In this section, we develop a method for deconvolution of the
blurred image with the translational and rotational motions. We
utilize the special property of the Laguerre-Fourier expansion in
SE(2) Fourier transform.

\subsubsection{Fourier transform of the Laguerre-Fourier expansion }
%%%%%%%% 3.3.1 %%%%%%%%%%%%%%%

When a function on polar coordinates, $\rho(r,\phi)$ is defined as
a truncated Laguerre-Fourier expansion as
$$
\rho(r,\phi) = \sum^{N}_{k=0}\sum^{k}_{l=-k}\tilde{\rho}_{kl}
\chi^*_{kl}(r,\phi),
$$
its Fourier transform in $\SE{2}$ is
$$
\hat{\rho}_{mn}(p) =  \int^{2 \pi}_{\theta=0}\int^{2 \pi}_{\phi=0}
\int^{\infty}_{r=0} \rho(r,\phi) i^{n-m} e^{i[m\theta +
(n-m)\phi]} J_{m-n}(p r) rdrd\phi d\theta
$$
$$
= \sum^{N}_{k=0}\sum^{k}_{l=-k}\tilde{\rho}_{kl} \int^{2
\pi}_{\theta=0}\int^{2 \pi}_{\phi=0} \int^{\infty}_{r=0} y_{kl}(r)
z^*_{l}(\phi) i^{n-m} e^{i[m\theta + (n-m)\phi]} J_{m-n}(p r)
rdrd\phi d\theta
$$
$$
= i^{n-m} \sum^{N}_{k=0}\sum^{k}_{l=-k}\tilde{\rho}_{kl}
\int^{\infty}_{r=0} y_{kl}(r) J_{m-n}(p r)rdr \int^{2
\pi}_{\phi=0} e^{il\phi} e^{i(n-m)\phi} d\phi \int^{2
\pi}_{\theta=0} e^{im\theta} d\theta
$$
$$
=4 \pi^2 i^{n-m} \sum^{N}_{k=0}\sum^{k}_{l=-k}\tilde{\rho}_{kl}
\left( \int^{\infty}_{r=0} y_{kl}(r) J_{m-n}(p r)rdr
\right)\delta_{l,-n} \delta_{m,0}
$$
$$
=4 \pi^2 i^{n} \sum^{N}_{k=0}\sum^{k}_{l=-k}\tilde{\rho}_{kl}
\left( \int^{\infty}_{r=0} y_{kl}(r) J_{l}(p r)rdr
\right)\delta_{l,-n} \delta_{m,0}
$$
On the other hand, we have a useful identity in \cite{cava_nume}
and\cite{book} as
$$
\int^{\infty}_{0} (\alpha r)^m L^m_n(\alpha^2 r^2)e^{-\alpha^2
r^2/2}J_m(kr)rdr = (-1)^n \alpha^{-2}(k/\alpha)^m
L^m_n(k^2/\alpha^2)e^{-k^2/2\alpha^2}.
$$
Using this identity, we can have
$$
\int^{\infty}_{r=0} y_{kl}(r) J_{l}(p r)rdr = (-1) ^
{(k-|l|)/2}\sqrt{\frac{[(k-|l|)/2]!}{\pi[(k+|l|)/2]!}}
\int^{\infty}_{r=0}    r^{|l|} L^{|l|}_{(k-|l|)/2}(r^2)e^{-r^2/2}
 J_{l}(p r)rdr
$$
$$
 = (-1) ^
{(k-|l|)/2}\sqrt{\frac{[(k-|l|)/2]!}{\pi[(k+|l|)/2]!}}
(-1)^{(k-l)/2} p^{|l|} L^{|l|}_{(k-|l|)/2} (p^2)e^{-p^2/2}
=(-1)^{(k-l)/2} y_{k,l}(p).
$$
Therefore,
\begin{equation}\label{eq:ft_rho_HL}
\hat{\rho}_{mn}(p) =4 \pi^2 i^{n}
\sum^{N}_{k=0}\sum^{k}_{l=-k}\tilde{\rho}_{kl}
(-1)^{{\frac{k-l}{2}}} y_{k,l}(p) \delta_{l,-n} \delta_{m,0} = 4
\pi^2 i^{n} \sum^{2[\frac{N-n}{2}]+n}_{k=|n|}\tilde{\rho}_{k,-n}
(-1)^{\frac{k+n}{2}} y_{k,-n}(p)
 \delta_{m,0},
\end{equation}
where $[n/2]=n/2$ if $n$ is even and $[n/2]=(n-1)/2$ if $n$ is
odd.

\subsubsection{Deconvolution using the Laguerre-Fourier expansion }
%%%%%%%% 3.3.2 %%%%%%%%%%%%%%%
%
%Using (\ref{eq:ft_conv_3}),(\ref{eq:ft_f3}) and
%(\ref{eq:ft_rho_HL}), we have
%$$
%\hat{\gamma}_{m,n}(p) = 4 \pi^2 i^{n}\delta_{m,0}
%\sum^{N}_{k=|n|}\tilde{\rho}_{k,-n} (-1)^{\frac{k+n}{2}}
%Y_{k,-n}(p)  e^{-p^2 t_1} e^{-n^2
% t_2}.
%$$
%Now, let us consider what $\gamma(r,\phi)$ looks like.

In Section \ref{sc:subsc_HL_case_1}, we noticed that the
translational blurred version of (\ref{eq:truncated_H}) is
(\ref{eq:blurred_H_case1}). This scaling effect on the domain
appears in the polar coordinates as follows: If the original image
is given as
\begin{equation} \label{eq:truncated_L_2}
\rho(r,\phi) =  \sum^{N}_{m=0} \sum^{m}_{n=-m}
 \tilde{\rho}_{mn}\chi^*_{mn}(r,\phi),
\end{equation}
then its translational motion blurring is
\begin{equation}\label{eq:blurred_L_case3}
\gamma (r,\phi) = \sum^{N}_{m=0} \sum^{m}_{n=-m}
\tilde{\gamma}_{mn}\chi^*_{mn}(\frac{r}{a},\phi),
\end{equation}
since the (\ref{eq:truncated_L_2}) and (\ref{eq:truncated_H}) are
equivalent under the simple coordinate relation, $r=\cos\phi$ and
$r=\sin\phi$\cite{park:interconversion}.

 In Section \ref{sc:subsc_HL_case_2}, the rotational blurred
version of (\ref{eq:truncated_L}) retains the structure of the
truncated Laguerre-Fourier expansion. In other words, if the
original image is expressed as a truncated Laguerre-Fourier
expansion, then its rotational blurring gives a truncated
Laguerre-Fourier expansion with the same truncation limit and
domain.

Therefore, we can conclude that the motion blurring in $\SE{2}$ of
(\ref{eq:truncated_L_2}) has the structure of
(\ref{eq:blurred_L_case3}), because the full motion in $\SE{2}$
can be decomposed into translation and rotation. The details of
the commutativity of the two motion blurring will be shown in the
appendix.

From (\ref{eq:ft_rho_HL}) and the convolution theorem, we have
$$
\hat{\gamma}_{m,n}(p) = \left(\hat{\rho}(p)
\hat{f_3}(p)\right)_{m,n} =  4 \pi^2 i^{n}\delta_{m,0}
\sum^{2[\frac{N-n}{2}]+n}_{k=|n|}\tilde{\rho}_{k,-n}
(-1)^{\frac{k+n}{2}} y_{k,-n}(p)  e^{-p^2 t_1} e^{-n^2
 t_2}.
$$
Also, we can have the Fourier transform of
(\ref{eq:blurred_L_case3}) directly as
$$
\hat{\gamma}_{m,n}(p) = 4 \pi^2 i^{n}\delta_{m,0}
\sum^{2[\frac{N-n}{2}]+n}_{k=|n|}\tilde{\gamma}_{k,-n}
(-1)^{\frac{k+n}{2}} a^2 y_{k,-n}(ap).
$$
Equating the two expressions of $\hat{\gamma}_{m,n}(p)$ gives
$$
\sum^{2[\frac{N-n}{2}]+n}_{k=|n|}\tilde{\rho}_{k,-n}
(-1)^{\frac{k}{2}} y_{k,n}(p) e^{-p^2 t_1} e^{-n^2
 t_2} =\sum^{2[\frac{N-n}{2}]+n}_{k=|n|}\tilde{\gamma}_{k,-n} (-1)^{\frac{k}{2}} a^2
y_{k,n}(ap)
$$
because $y_{m,n}=y_{m,-n}$. Since this should holds for any choice
of $p$ with a fixed $n$, we can have the following matrix
expression as in Section \ref{sc:subsc_HL_case_1}:
$$
WYJ \times r = a^2 Y_a J \times g,
$$
where $W_{k,l} = \delta_{k,l}e^{-(p^{(k)})^2 t_1}e^{-n^2 t_2}$,
 $J_{k,l} = \delta_{k,l} (-1)^{\frac{|n|}{2} + k -1}$,
$$
r = \left[\hspace{3mm} \tilde{\rho}_{|n|,-n} \hspace{3mm}
\tilde{\rho}_{|n|+2,-n} \hspace{3mm}\cdots \hspace{3mm}
\tilde{\rho}_{2[\frac{N-n}{2}]+n,-n} \hspace{3mm} \right]^T
$$
$$
g = \left[\hspace{3mm} \tilde{\gamma}_{|n|,-n} \hspace{3mm}
\tilde{\gamma}_{|n|+2,-n} \hspace{3mm}\cdots \hspace{3mm}
\tilde{\gamma}_{2[\frac{N-n}{2}]+n,-n} \hspace{3mm} \right]^T
$$
$$
Y=\left[
\begin{array}{ccc}
  y_{|n|,n}(p^{(1)}) & y_{|n|+2,n}(p^{(1)}) & \vdots \\
  y_{|n|,n}(p^{(2)}) & \ddots & \vdots \\
  \cdots & \cdots & y_{2[\frac{N-n}{2}]+n,n}(p^{(M)}) \\
\end{array}
\right]
$$
and
$$
Y_a=\left[
\begin{array}{ccc}
  y_{|n|,n}(ap^{(1)}) & y_{|n|+2,n}(ap^{(1)}) & \vdots \\
  y_{|n|,n}(ap^{(2)}) & \ddots & \vdots \\
  \cdots & \cdots & y_{2[\frac{N-n}{2}]+n,n}(ap^{(M)}) \\
\end{array}
\right]
$$
$J^{-1} = J$ and pseudo-inverse of $Y$ is given by $Y^{+}= (Y^T
Y)^{-1}Y^T$ if the sampling points, $p$'s are chosen appropriately
as shown in the previous work~\cite{park:interconversion}, which
guarantee the stable inversion of $Y$. Since inverting $W$ needs
regularization, $W^{+}_{k,l}=\delta_{k,l} ( 1 / (
e^{-t_1(p^{(k)})^2} e^{-t_2n^2 } + \epsilon )$ with a small
number, $\epsilon$.

Now we have
$$
r = a^2  J  Y^+ W^{+}  Y_a J \times g,
$$
Note that this holds for a fixed value for $n$. Therefore,
applying it for a different $n$ gives the full
$\tilde{\rho}_{m,n}$, which is the Laguerre-Fourier coefficients
for the deblurred image.

\section{Numerical Examples}\label{sc:numericalex}

\begin{figure}
  \centering
  \includegraphics[width=2.15in]{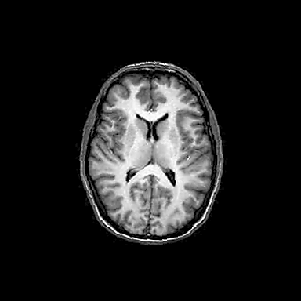}\\
  \caption{Original image. The original image is from http://www.med.univ-angers.fr/discipline /radiologie/Intlatlas/t1ax11.html}\label{fig:original}
\end{figure}

\begin{figure}
    \centering
        \subfigure[Blurred image with translational motion]{ \label{fig:bl_case_1}
            \includegraphics[width=2.15 in]{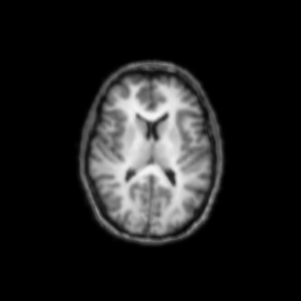}
             }
        \subfigure[Blurred image with rotational motion]{ \label{fig:bl_case_2}
            \includegraphics[width=2.15 in]{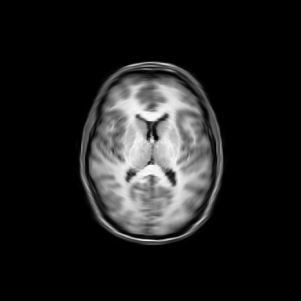}
            }
        \subfigure[Blurred image with translational and rotational motion]{ \label{fig:bl_case_3}
            \includegraphics[width=2.15 in]{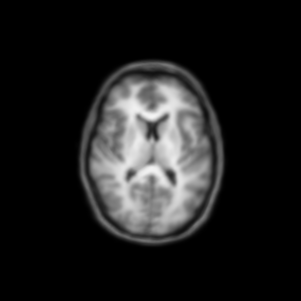}
            }
     \caption{Blurred image} \label{fig:blurred_images}
\end{figure}

We generate artificially blurred data as follows; We sample the
amount of motions from Gaussian motion distributions. Then the
original image is shifted by the sampled motion. We prepare a set
of the shifted images and average them to have our `experimental'
blurred image, $\gamma({\bf x})$. With knowledge of the motion
distribution $f(g)$, we examine if inverting will give back a good
estimate of the original. The original image and our artificially
blurred images are shown in Fig. \ref{fig:original} and Fig.
\ref{fig:blurred_images}, respectively. We averaged 100 shifted
images for each blurred image in Fig. \ref{fig:blurred_images}.

\subsection{Case 1 : Deblurring of translational blurring }
\label{sc:subsc_numerical_case_1}

The deconvolution method described in Section
\ref{sc:subsc_FT_case_1} can be implemented by discrete Fourier
transform(DFT). In order to capture the detail of the motion
distribution, $f_1$, we need to sample $f_1$ on a very fine grid,
because the distribution function is highly concentrated near the
mean value. To match the consistency DFTs of the samples of $f_1$
and the discrete image data, we have to resample the image on the
fine grid on which we sample $f_1$. This resampling is done by
simple linear interpolation of the image.

On the other hand, in order to implement the deconvolution method
using the Hermite expansion described in Section
\ref{sc:subsc_HL_case_1}, we have to have an appropriate truncated
Hermite expansion of the blurred image. The process to obtain a
truncated Hermite expansion of a discrete image was  developed in
\cite{park:interconversion}. Essentially, a truncated Hermite
expansion optimally fit to the image function is obtained. In
contrast to the Fourier method, this method gives the
interpolation values automatically, because the fitted truncated
expansion is a continuous function. Furthermore, in the
formulation, the blurred image function retains the structure of
the original expansion. After having the truncated Hermite
expansion, the estimation of the original image can be obtained by
the method in Section \ref{sc:subsc_HL_case_1}

\begin{figure}
    \centering
        \subfigure[Deconvolution of the translational blurred image using Fourier transform]{ \label{fig:de_bl_ft_case_1}
            \includegraphics[width=2.15 in]{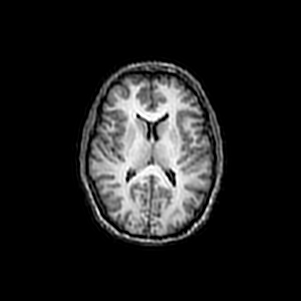}
             }
        \subfigure[Deconvolution of the translational blurred image using Hermite expansion]{ \label{fig:de_bl_hl_case_1}
            \includegraphics[width=2.15 in]{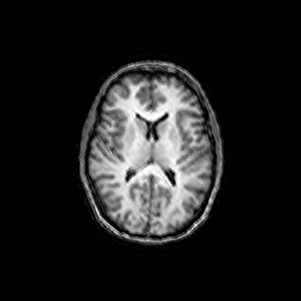}
             }
     \caption{Deconvolution of the translational blurred image } \label{fig:deblurred_images_case_1}
\end{figure}

Numerical results of the two methods for the deconvolution of the
translational blurring are shown in Fig
\ref{fig:deblurred_images_case_1}.

\subsection{Case 2 : Deblurring of rotational blurring }
\label{sc:subsc_numerical_case_2}

As the implementation of the Fourier method for deconvolution of
the translational motion blur in Section \ref{sc:subsc_FT_case_1},
the Fourier method deconvolution of the rotational motion blur in
Section \ref{sc:subsc_FT_case_2} also needs resampling of the
blurred image on a fine and polar grid, because the distribution
function, $f_2$ is highly concentrated and the formulation is done
on polar coordinates. Thus, interpolation of the blurred image on
a fine polar grid is performed to match the DFTs of the
distribution and the blurred image.

For the method using the Laguerre-Fourier expansion in Section
\ref{sc:subsc_HL_case_2}, firstly we have a truncated
Laguerre-Fourier expansion for the discrete blurred image.
Specifically a truncated Hermite expansion for the image is
obtained first and then we convert it to a truncated
Laguerre-Fourier expansion\cite{park:interconversion}. The
Laguerre-Fourier coefficients of the deblurred image can be
computed by (\ref{eq:hl_case_2_reg}).

\begin{figure}
    \centering
        \subfigure[Deconvolution of the rotational blurred image using Fourier transform]{ \label{fig:de_bl_ft_case_2}
            \includegraphics[width=2.15 in]{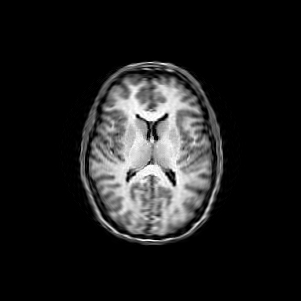}
             }
        \subfigure[Deconvolution of the rotational blurred image using Laguerre-Fourier expansion]{ \label{fig:de_bl_hl_case_2}
            \includegraphics[width=2.15 in]{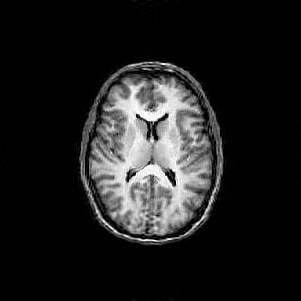}
             }
     \caption{Deconvolution of the rotational blurred image } \label{fig:deblurred_images_case_2}
\end{figure}

Numerical results of the two methods for the deconvolution of the
rotational blurring are shown in Fig
\ref{fig:deblurred_images_case_2}.

\subsection{Case 3 : Deblurring of combined translational and rotational blurring  }
\label{sc:subsc_numerical_case_3}

In order to have Fourier transform in $\SE{2}$ of the blurred
image, we resample it on fine polar grid. With this sampling, we
can compute the Fourier transform by the discretized version of
(\ref{eq:ft_se2}). The Fourier transform of the original image can
be estimated by (\ref{eq:ft_case_3_reg}).

As in the case 2, once we have a truncated Laguerre-Fourier
expansion for the blurred image, we can obtain the
Laguerre-Fourier coefficients of the deblurred image by the method
in Section \ref{sc:subsc_HL_case_3}, which is given as matrix
manipulations.

\begin{figure}
    \centering
        \subfigure[Deconvolution of  combined motional blurred image using Fourier transform]{ \label{fig:de_bl_ft_case_3}
            \includegraphics[width=2.15 in]{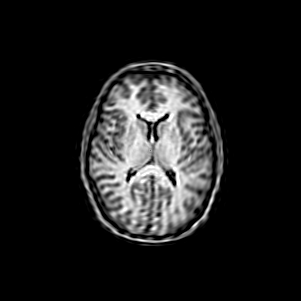}
             }
        \subfigure[Deconvolution of  combined motional blurred image using Laguerre-Fourier expansion]{ \label{fig:de_bl_hl_case_3}
            \includegraphics[width=2.15 in]{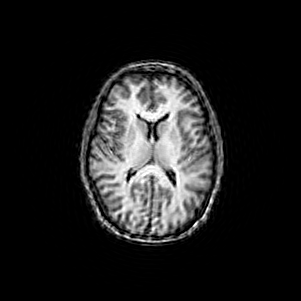}
             }
     \caption{Deconvolution of combined translational and rotational blurred image } \label{fig:deblurred_images_case_3}
\end{figure}

Numerical results of the two methods for the deconvolution of the
translational and rotational blurring are shown in Fig
\ref{fig:deblurred_images_case_3}.

%\section{Conclusions}
%
%Methods for solving the motional deblurring problem are developed
%here.

\section{Conclusion and discussion}
In this work, we have shown our method to restore the 2D blurred
images that were generated by the translational motion, the
rotational motion, or both. In our formulation, the blurring
process can be expressed as a convolution of an original image and
the motion distribution function. Since the convolution can be
interpreted as the multiplication in the Fourier space, the
deconvolution ( i.e. deblurring ) is the simple inversion process
in that space.

We applied this concept to the three cases. We used the Fourier
analysis in $\SE2$ for the combined translational and rotational
motion blur, which is defined using the representation theory of
the Euclidean motion group of the plane.

The Fourier analysis has the strong advantage that the convolution
in the original space can be replaced with the multiplication in
the Fourier space, which is much simpler to manipulate. However
the implementation of it needs more manipulation of data, because
the image is the function defined only on a discrete grid domain,
while the formulation is done in the continuous domain. Although
the discrete Fourier transform is a good way to do it, resampling
is required before the transform as described in Section
\ref{sc:numericalex}. Shortly, the resampling is needed for
matching the resolution of the image function and the distribution
function. In addition to the resolution matching, we need another
resampling process for the cases of the rotational motion blur and
the combined motion blur, since the image function on the polar
grid is required in these cases.

To overcome this, we proposed an alternative method using the
Hermite and Laguerre-Fourier expansions. Once the 2D truncated
Hermite expansion optimally fitted to the given discrete image is
obtained, we can directly utilize the results of the derivation
without introducing discrete version of it, which is essential in
the Fourier method. In this case, resampling on a finer grid comes
naturally, because our expansion is already a function on a
continuous domain. On the other hand, we already know that the
interconversion of the Hermite and Laguerre-Fourier expansions is
possible losslessly and it can be viewed as the coordinate
conversion from Cartesian to polar
coordinates~\cite{park:interconversion}. Therefore, once we have
the Hermite expansion for an image, the corresponding expression
on the polar coordinates, which is exactly the Laguerre-Fourier
expansion, can be obtained easily. Since the two expansions have
the special property under the Fourier transform, which is that
they retain their structure under the Fourier transform, the
expansions enable the straightforward implementation of the
deconvolution method.

The computational complexity of the interconversion of the two
coordinates is still high($\mathcal{O}(n^4)$). However, the
aforementioned advantages of the Hermite and Laguerre-Fourier
expansions are still valid and we leave the improvement of the
conversion algorithm for the future work.

For both of the two approaches for deconvolution, we need some
objective measurement scheme to assess the deblurred data in order
to pick the best regularization parameter that gives the best
deblurred image. With various values of the regularization
parameter, we have a set of candidates for the estimation of the
original image. In our work, we choose one with naked eye, which
is not too blurry or too sharp. For a more systematic and
objective assessment, a special kind of measuring method should be
proposed. Furthermore, if we don't know the variance of the
Gaussian motion distribution, which is a more practical situation,
our tunable parameters are two: variance and regularization
parameter. We also leave the development of an objective tuning
algorithm for the future.

\section{Appendix}%

 On a polar plane and a circle, the Gaussian distribution
densities are
$$ f_1(r,\phi) = \frac{1}{4 \pi t_1}e^{-r^2/4t_1}
\hspace{5mm} \mbox{  and  } \hspace{5mm} f_2(\theta) = \frac{1}{2
\pi }\sum^{\infty}_{k=-\infty}e^{-k^2 t_2}e^{ik\theta}
$$
Similarly we define the functions in $\SE2 $ as
$$
F_1(r,\phi,\theta) = \frac{1}{4 \pi
t_1}e^{-r^2/4t_1}\delta(\theta) \hspace{5mm} \mbox{  and  }
\hspace{5mm} F_2(r,\phi,\theta) = \frac{1}{2 \pi
}\sum^{\infty}_{k=-\infty}e^{-k^2 t_2}e^{ik\theta}
\frac{\delta(r)}{2\pi r} \delta(\phi).
$$
We want to show that
$$
(F_1*F_2)(g) = (F_2*F_1)(g) = f_1(g) f_2(g).
$$

\vspace{7mm}
\noindent{\bf Proof}\\
\noindent
$$
(F_1*F_2)(g) = \int_{\SE2} F_1(h)F_2(h^{-1}\circ g)dh
$$
and
$$
(F_2*F_1)(g) = \int_{\SE2}  F_2(h)F_1(h^{-1}\circ g)dh.
$$
By changing the variable $k=h^{-1}\circ g$, we have
$$
(F_1*F_2)(g) = \int_{\SE2}  F_1(g\circ k^{-1} )F_2(k)dk=
\int_{\SE2} F_2(h)F_1(g\circ h^{-1} )dh.
$$

\noindent $g$ and $h$ can be parameterized as
$$
g=g(r,\phi,\theta) =\left[
\begin{array}{ccc}
    \cos \theta & -\sin \theta & r\cos \phi \\
    \sin \theta  & \cos \theta & r\sin \phi \\
    0 & 0 & 1 \\
  \end{array}\right]
\hspace{6mm} \mbox{and} \hspace{6mm}
 h=h(R,\Phi,\Theta)
=\left[
\begin{array}{ccc}
    \cos \Theta & -\sin \Theta & R\cos \Phi \\
    \sin \Theta  & \cos \Theta & R\sin \Phi \\
    0 & 0 & 1 \\
  \end{array}\right]
$$
The multiplication is
$$
h^{-1}\circ g =\left[
\begin{array}{ccc}
    \cos (\theta-\Theta) & -\sin (\theta-\Theta)  & r\cos (\phi-\Theta) -R \cos (\Phi-\Theta)  \\
    \sin (\theta-\Theta) & \cos  (\theta-\Theta) & r\sin (\phi-\Theta) -R \sin (\Phi-\Theta))\\
    0 & 0 & 1 \\
  \end{array}\right]
$$
$$
g\circ   h^{-1}=\left[
\begin{array}{ccc}
    \cos (\theta-\Theta) & -\sin (\theta-\Theta)  & r\cos  \phi  -R \cos (\theta+\Phi-\Theta)  \\
    \sin (\theta-\Theta) & \cos  (\theta-\Theta) & r\sin  \phi -R \sin (\theta+\Phi-\Theta)\\
    0 & 0 & 1 \\
  \end{array}\right]
$$
Therefore,
$$
(F_1*F_2)(g) = \int F_2(h)F_1(g\circ h^{-1} )dh
$$
$$
= \int_{\SE2}  \left\{\frac{1}{2 \pi
}\sum^{\infty}_{k=-\infty}e^{-k^2 t_2}e^{ik\Theta}
\frac{\delta(R)}{2\pi R} \delta(\Phi) \right\}\left\{  \frac{1}{4
\pi t_1}e^{-( R^2+r^2-2Rr\cos (\phi-\theta-\Phi+\Theta) )
/4t_1}\delta(\theta-\Theta) \right\} R \hspace{0.1mm}dR
\hspace{0.1mm} d\Phi \hspace{0.1mm}d\Theta.
$$
Integration over $\Theta$ gives
$$
 = \int
\left\{\frac{1}{2 \pi }\sum^{\infty}_{k=-\infty}e^{-k^2
t_2}e^{ik\theta} \frac{\delta(R)}{2\pi R}  \delta(\Phi)
\right\}\left\{ \frac{1}{4 \pi t_1}e^{-( R^2+r^2-2Rr\cos (\phi
-\Phi ) ) /4t_1} \right\} R \hspace{0.1mm}dR \hspace{0.1mm} d\Phi
\hspace{0.1mm}
$$
$$
 =
 \frac{1}{8 \pi^2 t_1}\sum^{\infty}_{k=-\infty}e^{-k^2
t_2}e^{ik\theta}      \int \left\{ e^{-( R^2+r^2-2Rr\cos (\phi
-\Phi ) ) /4t_1} \frac{\delta(R)}{2\pi R}  \delta(\Phi) \right\} R
\hspace{0.1mm}dR \hspace{0.1mm} d\Phi
$$
Using the fact that the $ \delta(R)/(2\pi R)$ is a special delta
function on a polar coordinate at singularity ( $R=0$ ), we have
$$
 (F_1*F_2)(g)=
 \frac{1}{8 \pi^2 t_1}\sum^{\infty}_{k=-\infty}e^{-k^2
t_2}e^{ik\theta}       e^{-  r^2  /4t_1}
$$

Similarly,
$$
(F_2*F_1)(g) = \int F_2(h)F_1(h^{-1}\circ g)dh
$$
$$
= \int_{\SE2}  \left\{\frac{1}{2 \pi
}\sum^{\infty}_{k=-\infty}e^{-k^2 t_2}e^{ik\Theta}
\frac{\delta(R)}{2\pi R} \delta(\Phi)  \right\}\left\{  \frac{1}{4
\pi t_1}e^{-( R^2+r^2-2Rr\cos (\phi -\Phi ) )
/4t_1}\delta(\theta-\Theta) \right\} R \hspace{0.1mm}dR
\hspace{0.1mm} d\Phi \hspace{0.1mm}d\Theta.
$$
$$
=
 \frac{1}{8 \pi^2 t_1}\sum^{\infty}_{k=-\infty}e^{-k^2
t_2}e^{ik\theta}       e^{-  r^2  /4t_1}
$$

Therefore, we showed that
$$
(F_1*F_2)(g) = (F_2*F_1)(g) = f_1(g) f_2(g)=\frac{1}{8 \pi^2
t_1}\sum^{\infty}_{k=-\infty}e^{-k^2 t_2}e^{ik\theta}       e^{-
r^2  /4t_1},
$$
where
$$
g=g(r,\phi,\theta) =\left[
\begin{array}{ccc}
    \cos \theta & -\sin \theta & r\cos \phi \\
    \sin \theta  & \cos \theta & r\sin \phi \\
    0 & 0 & 1 \\
  \end{array}\right]
  $$

\end{document}